\documentclass[12pt]{article} 
\usepackage{amsfonts,amsmath,layout}
\usepackage{graphicx}   
\usepackage{url}
\usepackage{hyperref}

\newcommand{\be}{\begin{equation}}
\newcommand{\ee}{\end{equation}}

\newcommand{\N}{\mathbb N}

\newcommand{\R}{\mathbb R}

\def\QED{\hfill$\; \Box$\medskip}

\def\R{\mathbb R}
\def\N{\mathbb N}

\def\E{\mathbb E}
\def\P{\mathbb P}

\def\Om{\Omega}
\def\lg{\langle}
\def\rg{\rangle}
\def\F{\mathcal F}
\def\FR{{\cal F}}

\def\BN{{\bf B}}
\def\ep{\epsilon}

\def\ds{\displaystyle}

\newtheorem{Theorem}{Theorem}[section]

\newtheorem{Proposition}[Theorem]{Proposition}
\newtheorem{Lemma}[Theorem]{Lemma}

\newtheorem{Remark}[Theorem]{Remark}

\title{{H\"{o}lder regularity for viscosity solutions of fully nonlinear, local or nonlocal, Hamilton-Jacobi equations
with super-quadratic growth in the gradient}
\thanks{This work was partially supported by the French ANR (Agence Nationale de la Recherche) through MICA project (ANR-06-06-BLAN-0082) and KAMFAIBLE project (ANR-07-BLAN-0361).}}
\author{Pierre Cardaliaguet and Catherine Rainer \thanks{Universit\'e de Bretagne
Occidentale, UMR 6205, 6 Av. Le Gorgeu,
BP 809, 29285 Brest (France); e-mail:
[Pierre.Cardaliaguet, Catherine.Rainer]@univ-brest.fr}}

\begin{document}
\maketitle

\begin{abstract} Viscosity solutions of fully nonlinear, local or non local, Hamilton-Jacobi equations with a super-quadratic growth 
in the gradient variable are proved to be H\"{o}lder continuous, with a modulus depending only on the growth of the Hamiltonian. The proof involves
some representation formula for nonlocal Hamilton-Jacobi equations in terms of controlled jump processes and a weak reverse inequality. 
\end{abstract} 

\medskip
\noindent
{\bf \underline{Key words}:} Integro-Differential Hamilton-Jacobi equations, viscosity solutions, H\"older continuity, 
degenerate parabolic equations, reverse H\"older inequalities.

\medskip
\noindent
{\bf \underline{MSC Subject classifications}:}  49L25, 35K55, 93E20, 26D15  \\


\section{Introduction}

In a previous paper \cite{CC}, the first author investigated the regularity of solutions to the Hamilton-Jacobi equation
\be\label{HJ0}
u_t(x,t)-{\rm Tr}\left(a(x,t)D^2u(x,t)\right)+H(x,t,Du(x,t))=0\qquad {\rm in }\; \R^N\times (0,T)
\ee
under a super-quadratic growth condition on the Hamiltonian $H$ with respect to the gradient variable: 
$$
\frac{1}{\delta}|z|^q-\delta \leq H(x,t,z) \leq \delta|z|^q+\delta\qquad \forall (x,t,z)\in \R^N\times(0,T)\times \R^N\;,
$$
for some $\delta\geq 1$, $q\geq 2$. Under this assumption, it is proved in \cite{CC} that
any continuous, bounded solution $u$ of \eqref{HJ0} is H\"{o}lder continuous
on $\R^N\times [\tau,T]$ (for any $\tau\in (0,T)$), with H\"{o}lder exponent and constant only 
depending on $N$, $p$, $\delta$,
 $q$, $\tau$ and  $\|u\|_\infty$. In particular, such a modulus of continuity is independent of the regularity of $a$ and of $H$ with respect
 to the variables $(x,t)$. The result is somewhat surprizing since no uniform ellipticity on $a$ is required. \\

The aim of this paper is to extend this regularity result to solutions of fully nonlinear, local or nonlocal, Hamilton-Jacobi equations
which have a super-quadratic growth with respect to the gradient variable. Beside its own interest, 
such a uniform estimate is important in homogenization theory where, for instance, it is used to prove the existence of correctors. \\
 
Let us consider a fully nonlinear, nonlocal Hamilton-Jacobi equation of the form
\be\label{HJNL0}
u_t+F(x,t,Du,[u])=0\qquad {\rm in }\; \R^N\times (0,T)\;.
\ee
In the above equation we assume that the mapping $F:\R^N\times (0,T)\times \R^N\times {\mathcal C}^2_b(\R^N)\to \R$ is 
nonincreasing with respect to the nonlocal variable,
i.e., 
$$
\left[\phi\leq \psi \; {\rm and }\; \phi(x)=\psi(x) \right]\; \Rightarrow \; F(x,t,\xi,[\phi]) \geq F(x,t,\xi,[\psi])
$$
for any function $\phi,\psi\in {\mathcal C}^2_b(\R^N)$. 
Let us recall (see \cite{AT,Ar,BBP,BI, Ph,Sa1,Sa2, So} for instance)
 that a subsolution (resp. a supersolution) of equation (\ref{HJNL0}) is a continuous map $u:\R^N\times[0,T]\to \R$ such that, for  
any continuous, bounded test function $\phi:\R^N\times (0,T)\to\R$, which has continuous second order derivatives and such that 
$u-\phi$ has a global maximum (resp. global minimum)  at some point $(\bar x,\bar t)$, one has 
$$
\phi_t(\bar x,\bar t)+F(x,t,D\phi(\bar x,\bar t),[\phi(\cdot,\bar t)])\leq 0\qquad ({\rm resp.}\; \geq 0) \;.
$$

Our main assumption on $F$ is the following structure condition, which roughly says that 
$F$ is super-quadratic with respect to the gradient variable:
\be\label{StrucCondNL}
-\delta M^+[\phi](x) +\frac{1}{\delta}|\xi|^q -\delta \ \leq \
F(x,t,\xi,[\phi])  \ \leq  \
-\delta  M^-[\phi](x) +\delta|\xi|^q +\delta
\ee
for any $(x,t,\xi,\phi)\in \R^N\times (0,T)\times \R^N\times {\mathcal C}^2_b(\R^N)$, for some constants $q>2$ and 
$\delta\geq 1$, where $M^-$ and $M^+$ are defined by
$$
M^-[\phi](x)=\inf_{\lambda\in (0,1],\ b\in \BN\backslash \{0\} } \left\{ \frac{\phi(x+\lambda b)-\phi(x)-\lg D\phi(x), \lambda b\rg}{|b|^2} \right\}
$$
and
$$
M^+[\phi](x)=\sup_{\lambda\in (0,1],\ b\in \BN\backslash \{0\}} \left\{ \frac{\phi(x+\lambda b)-\phi(x)-\lg D\phi(x), \lambda b\rg}{|b|^2} \right\}
$$
and where $\BN$ is the unit ball of $\R^N$. Let us note that, under the above assumption, a solution of (\ref{HJNL0}) is a supersolution of 
\be\label{intro:HJnl+}
u_t-\delta M^-[u(\cdot,t)](x) +\delta|Du|^q +\delta=0\qquad {\rm in }\; \R^N\times (0,T)
\ee
and a subsolution of 
\be\label{intro:HJnl-}
u_t-\delta M^+[u(\cdot,t](x) +\frac{1}{\delta}|Du|^q -\delta=0\qquad {\rm in }\; \R^N\times (0,T)\;.
\ee
We denote by $p$  the conjugate exponent of $q$. 
We are interested in solutions which are bounded by some constant $M$:
\be\label{DefM}
|u(x,t)|\leq M \qquad \forall (x,t)\in \R^N\times [0,T]\;.
\ee
In what follows, a {\em(universal) constant} is a positive number depending on the given data $q,\delta, N$ 
and $M$ only. Universal constants will be typically labeled with  $C$, but also with different letters (e.g., $\theta, A,\dots$). 
Dependence on extra quantities will be accounted for by using parentheses (e.g., $C(r)$ denotes a constant depending also on $r$). The constant $C$ appearing in the proofs may change from line to line.

\begin{Theorem} \label{Regu2NL} Let $u\in {\mathcal C}(\R^N\times[0,T])$ be a viscosity supersolution of \eqref{intro:HJnl+} and a subsolution
of (\ref{intro:HJnl-}), such that 
$|u|\leq M$ in $\R^N\times[0,T]$. For any $\tau\in (0,T)$, there are constants $\theta=\theta(\delta,M, N,q)>p$ and $C(\tau)= (\tau,\delta,M,N,q)>0$ 
such that
\begin{equation}\label{eq:main3NL}
|u(x_1,t_1)-u(x_2,t_2)|\leq C(\tau)\left[|x_1-x_2|^{(\theta-p)/(\theta-1)}+|t_1-t_2|^{(\theta-p)/\theta}\right]
\end{equation}
for any $(x_1,t_1), (x_2,t_2)\in \R^N\times [\tau, T]$.
\end{Theorem}

The main point of the above result is that \eqref{eq:main3NL} holds true uniformly with respect to 
$F$, as long as conditions \eqref{StrucCond} and the bound $|u|\leq M$ are  satisfied. 
In particular,  $\theta$ and  $C(\tau)$ are independent of the continuity modulus of $F$. In contrast to \cite{CS, Si}, $F$ can also be degenerate parabolic.\\

Let us note that the above result also applies to the solutions of the fully nonlinear, local equation
\be\label{HJloc}
u_t+F(x,t,Du,D^2u)=0\qquad {\rm in }\; \R^N\times (0,T)
\ee
provided $F:\R^N\times (0,T)\times \R^N\times {\mathcal S}^N\to \R$ is nonincreasing with respect to the 
matrix variable and satisfies the following structure condition:
\be\label{StrucCond}
-\delta\Lambda(X) +\frac{1}{\delta}|\xi|^q -\delta \leq 
F(x,t,\xi,X) \leq 
-\delta\lambda(X)+\delta|\xi|^q +\delta
\ee
for any $(x,t,\xi,X)\in \R^N\times (0,T)\times \R^N\times {\mathcal S}^N$, for some constants $q>2$ and 
$\delta\geq 1$, where
$$
\Lambda(X)= \max_{|z|\leq 1}\lg Xz,z\rg \;{\rm and }\; \lambda(X)=\min_{|z|\leq 1}\lg Xz,z\rg\qquad \forall X\in {\mathcal S}^N\;,
$$
(${\mathcal S}^N$ being the set of $N\times N$ symmetric matrices). Indeed, since
$$
M^-[\phi](x) \leq \lambda(D^2\phi(x))\qquad {\rm and }\qquad M^+[\phi](x) \geq \Lambda(D^2\phi(x))\;,
$$
any solution of (\ref{HJloc}) is a supersolution of (\ref{intro:HJnl+}) and a subsolution of (\ref{intro:HJnl-}). 
Note that 
$F$ is neither  required to be concave nor convex with respect to the matrix variable. \\

Here are some examples of nonlinear, nonlocal Hamilton-Jacobi equations satisfying the structure condition (\ref{StrucCondNL})
(see \cite{BBP} or \cite{BI} for instance): let us assume that 
$$
F(x,t,\xi,[\phi])= {\mathcal I}[\phi](x,t)+H(x,t,\xi)
$$
where $H$ is a  first order term with superquadratic growth:
$$
\frac{1}{\delta}|\xi|^q -\delta \leq H(x,t,\xi) \leq \delta|\xi|^q +\delta
$$
and where the nonlocal term ${\mathcal I}$ can be of the form
$$
{\mathcal I}[\phi](x,t)=\inf_{\alpha\in A}\sup_{\beta\in B}
\int_{\R^N} \phi(x+j_{\alpha,\beta}(x,t,e))-\phi(x)-\lg D\phi(x),j_{\alpha,\beta}(x,t,e)\rg \ d\nu (e)
$$
where $A, B$ are some sets, $j_{\alpha,\beta}:\R^N\times (0,T)\times \R^N\to \R^N$ is such that
$$
|j_{\alpha,\beta}(x,t,e)|\leq C(|e|\wedge 1) \qquad \forall (x,t,e,\alpha,\beta)\in \R^N\times (0,T)\times \R^N\times A\times B
$$
and the measure $\nu$ satisfies
\be \label{Intlambda}
\int_{\R^N} |e|^2\wedge 1\ d\nu(e) \leq C\;,
\ee
or  of the form
$$
{\mathcal I}[\phi](x,t)=\inf_{\alpha\in A}\sup_{\beta\in B}
\int_{\R^N} \phi(x+j_{\alpha,\beta}(x,t,e))-\phi(x)-\lg D\phi(x),j_{\alpha,\beta}(x,t,e)\rg {\bf 1}_{\BN}(e)\ d\nu (e)
$$
where $j_{\alpha,\beta}:\R^N\times (0,T)\times \R^N\to \R^N$ is now of linear growth
$$
|j_{\alpha,\beta}(x,t,e)|\leq C|e| \qquad \forall (x,t,e,\alpha,\beta)\in \R^N\times (0,T)\times \R^N\times A\times B
$$
and the measure $\nu$ again satisfies the integrability condition (\ref{Intlambda}). In this later case,
the part
$$
\int_{\BN} \phi(x+j_{\alpha,\beta}(x,t,e))-\phi(x)-\lg D\phi(x),j_{\alpha,\beta}(x,t,e)\rg {\bf 1}_{\BN}(e)\ d\nu (e)
$$
can be estimated from above and below by $M^+[\phi]$ and $M^-[\phi]$, while the part
$$
\int_{\R^N\backslash \BN} \phi(x+j_{\alpha,\beta}(x,t,e))-\phi(x)\ d\nu (e)
$$
is can be bounded---in equation (\ref{HJNL0})---by $CM$, where $M=\|u\|_\infty$ and $C$ is universal. \\

Some comments on the proof of Theorm \ref{Regu2NL} are now in order. 
As in \cite{CC}, the main ingredients are representation formulae for
simplified Hamilton-Jacobi equations, existence of ``nearly optimal trajectories", use of Brownian bridges and, 
finally, application of a reverse H\"{o}lder inequality. However, since we work with fully nonlinear, nonlocal equations, 
each step is technically more involved: the representation formulae (see Proposition \ref{prop:repHJnl+} or the proof of Proposition 
\ref{prop:SuperSol-nl}) are inspired by a work on controlled structure equations by the second author \cite{BMR}. They
involve controlled jump processes in a particular form. The estimates of the subsolutions in Proposition \ref{prop:SuperSol-nl}, 
which, in \cite{CC}, are obtained by controls
issued from Brownian bridge techniques, have to be built here in a much more subtle way: 
indeed the Hamiltonian of equation (\ref{intro:HJnl+}) being non convex, 
the naturally associated control problem should actually be a differential game, which has never been investigated in this framework. 
We overcome this difficulty by building
 explicit feedbacks. Finally, the construction of optimal trajectories in Lemma \ref{lem:IneqStochnl}, 
requires careful estimates because we are dealing with jump processes.\\

{\bf Notations : } Throughout the paper, $\BN$ denotes the closed unit ball of 
$\R^N$, ${\mathcal B}(\BN\backslash \{0\})$ the set of Borel measurable subsets of $\BN\backslash \{0\}$, 
${\mathcal C}(\R^N\times [0,T])$ the set of continuous functions on $\R^N\times [0,T]$ and
${\mathcal C}^2_b(\R^N)$ the set of bounded continuous functions on $\R^N$ with continuous second order  derivatives.

\section{Analysis of supersolutions}

Let $u$ be as in Theorem \ref{Regu2NL}. 
Throughout the proof of Theorem \ref{Regu2NL} it will be more convenient to work with $u(x,T-t)$ instead of $u(x,t)$. We  note that 
$u(x,T-t)$ is a supersolution of 
\be\label{HJnl+}
-v_t-\delta M^-[v(\cdot,t)](x) +\delta|Dv|^q +\delta=0\qquad {\rm in }\; \R^N\times (0,T)
\ee 
and a subsolution of
\be\label{HJnl-}
-v_t-\delta M^+[v(\cdot,t)](x) +\frac{1}{\delta}|Dv|^q -\delta=0\qquad {\rm in }\; \R^N\times (0,T)\;.
\ee
To simplify the notation we will write $u(x,t)$ instead of $u(x,T-t)$. \\

In this part we are concerned with some monotonity property along particular trajectories of supersolution of equation
(\ref{HJnl+}). For this we have to give a representation formula for solutions of this equation in terms of controlled jump processes. \\ 

Let $(\Om, \F,\P)$ be a complete probability space on which  is defined 
a $N-$dimensional Poisson random measure $\mu$. We assume that the Levy measure of $\mu$, denoted by 
$\nu$, is supported in $\BN$, has no atom and satisfies the conditions
\begin{equation}
\label{poisson}
\int_{\BN} |e|^2 \; d\nu(e) <+\infty \qquad{\rm and }\qquad \nu({\bf B})=+\infty\;.
\end{equation}
We denote by $\tilde \mu(de,dt)= \mu(de,dt)-\nu(de) dt$ the compensated Poisson measure. For all $t\in[0,T]$, $(\FR_{t,s},s\in [t,T])$ will denote the filtration generated by $\mu$ on the interval $[t,T]$, i.e. $\FR_{t,s}=\sigma\{ \mu([t,r]\times A),t\leq r\leq s, A\in{\cal B}(\R^N)\}$ completed by all null sets of $P$.\\

Let ${\mathcal A}(t)$ be the set of $(\FR_{t,s})$-adapted controls $(a_s)=(\lambda_s, b_s):[t,T]\to (0,1]\times {\bf B}\backslash\{0\}$
and let $L^p_{\rm ad}([t,T])$ be the set of $(\FR_{t,s})$-adapted controls $(\zeta_s):[t,T]\to \R^N$ such that
$\E\left[ \int_{t}^T |\zeta_s|^pds\right]<+\infty$. \\
To all control $(a_s)=(\lambda_s, b_s)\in {\mathcal A}(t)$ we associate a $(\FR_{t,s})$-martingale $M^a$ in the following way:\\
First we set
$$
\rho_s=\inf\left\{r>0\;;\; \nu(\BN\backslash B(0,r)) \leq \delta\frac{1}{|b_s|^2}\right\}\;,\; s\in[t,T].
$$
By the assumptions (\ref{poisson}) on the measure $\nu$, the process $(\rho_s)$ is well defined. It  takes its values in $[0,1]$ and is adapted to the filtration $(\F_{t,s})$.\\
Then we introduce  $A_s=A_s(a)=\BN\backslash B(0,\rho_s)$.
For all $s\in[0,T]$, the set $A_s$ belongs to ${\mathcal B}({\bf B}\backslash\{0\})\otimes \F_{t,s}$ and, since $\nu$ has no atoms, it satisfies
\be\label{As}
\nu(A_s)= \frac{\delta}{|b_s|^2}\qquad \mbox{\rm for almost all $s\in [t_,T]$, $\P-$a.s.}
\ee
We finally denote by $M^a$ the controlled martingale
\be\label{DefMa}
M^a_s=\int_{t}^s \int_{{\bf B}} \lambda_r b_r {\bf 1}_{A_r}(e) \tilde \mu(de,dr)\;.
\ee
Let us precise It\^{o}'s formula satisfied by $M^a$: for any smooth function $\phi\in {\mathcal C}^2_b$, 
$$
\begin{array}{l}
\ds{  \E[\phi(M^a_s)] }\\
\quad {\ds = \phi(0)+\E\left[ \int_{t}^s \int_{{\bf B}} \left(\phi(M^a_{r^-}+\lambda_r b_r{\bf 1}_{A_r}(e))-\phi(M^a_{r^-})
-\lg D\phi(M^a_{r^-}), \lambda_r b_r{\bf 1}_{A_r}(e)\rg \right)d\nu(e)dr \right]  }\\
\quad \ds{ = \phi(0)+\delta\E\left[ \int_{t}^s (\phi(M^a_{r^-}+\lambda_r b_r)-\phi(M^a_{r^-})
-\lg D\phi(M^a_{r^-}), \lambda_r b_r\rg)\frac{dr}{|b_r|^2}  \right] }
\end{array}
$$
Now we consider the  controlled system
\be\label{ContSystHJnl+}
\left\{\begin{array}{l}
dY_s= \zeta_sds+ dM^a_s\; , \; s\in[t,T],\\
Y_t=x\; ,
\end{array}\right.
\ee
where $\zeta\in L^p_{\rm ad}([t,T])$ and $a\in  {\mathcal A}(t)$.
The system (\ref{ContSystHJnl+}) is related to equation (\ref{HJnl+}) by the following proposition:

\begin{Proposition}\label{prop:repHJnl+} Let $v:\R^N\times [0,T]$ be a continuous viscosity solution to (\ref{HJnl+}). Then
$$
v(x,t)= \inf_{(\zeta,a)\in L^p_{\rm ad}([t,T])\times   {\mathcal A}(t)}
\E\left[ \ v(Y_T^{x,t,\zeta,a},T)+ C_{_+}\int_t^T |\zeta_s|^pds -\delta(T-t)\right]
$$
where $(Y_s^{x,t,\zeta,a})$ is the solution to (\ref{ContSystHJnl+}) and $C_{_+}>0$ is the  universal constant given by 
\begin{equation}\label{eq:C+}
C_{_+}=\frac{\delta^{-p/q}}{pq^{p/q}}\,.
\end{equation}
 \end{Proposition}

{\bf Proof : } The proof relies on two arguments: first the map
$$
w(x,t)=\inf_{(\zeta,a)\in L^p_{\rm ad}([t,T])\times  {\mathcal A}(t)}
\E\left[ \ v(Y_T^{x,t,\zeta,a},T)+ C_{_+}\int_t^T |\zeta_s|^pds -\delta(T-t)\right]
$$
is a viscosity solution of (\ref{HJnl+}). This result has been proved---in a slightly different framework---in \cite{BMR}, 
and we omit the proof, which is very close to that of \cite{BMR}. 
Second, in order to conclude that $w=v$, we need a uniqueness argument for the solution of (\ref{HJnl+}) with terminal condition
$v(\cdot,T)$. This is a direct consequence of the following comparison principle.
\hfill \QED

\begin{Lemma}[Comparison]\label{compar} Let $u$ be a continuous, bounded subsolution of (\ref{HJnl+}) on $\R^N\times [0,T]$ and
$v$ be a continuous, bounded supersolution of (\ref{HJnl+}). If $u(\cdot,T)\leq v(\cdot,T)$, then 
$u\leq v$ on $\R^N\times [0,T]$.
\end{Lemma}

{\bf Proof : } We use ideas of \cite{BI} for the treatment of the nonlocal term and of \cite{dll}
for the treatment of the super-linear growth with respect to
the gradient variable. Let $M=\max\{\|u\|_\infty, \|v\|_\infty\}$. Our aim is to show that, for any $\mu\in(0,1)$, 
\be\label{muuv}
\tilde u(x,t):= \mu u(x,t) -(1-\mu) M -(1-\mu)\delta(T-t)\; \leq\;  v(x,t)\;{\rm in }\; \R^N\times [0,T]\;.
\ee
Note that this inequality holds at $t=T$. One also easily checks that $\tilde u$ is a subsolution of equation
\be\label{HJnl+bis}
-w_t-\delta M^-[w(\cdot,t)](x) +\delta\mu^{1-q}|Dw|^q +\delta=0\qquad {\rm in }\; \R^N\times (0,T)\;.
\ee
For any $\ep>0$, let $\tilde u^\ep$ be the space-time sup-convolution of $\tilde u$ and $v_\ep$ be the  space-time inf-convolution of $v$ (see \cite{CIL}):
$$
 \tilde u^\ep(x,t)= \sup_{(y,s)\in \R^N\times [0,T]} \left\{ \tilde u(y,s)-\frac{1}{\epsilon}|(x,t)-(y,s)|^2 \right\}
$$
and
$$
v_\ep(x,t)= \inf_{(y,s)\in \R^N\times [0,T]} \left\{ v(y,s)+\frac{1}{\epsilon}|(x,t)-(y,s)|^2 \right\}\;.
$$
It is known that $\tilde u^\ep$ is still a subsolution of (\ref{HJnl+bis}) 
in $\R^N\times ((2M\ep)^\frac12, T-(2M\ep)^\frac12)$ while $v_\ep$ a supersolution of (\ref{HJnl+}) and that $u^\ep$ is semiconvex while
$v^\ep$ is semiconcave $\R^N\times ((2M\ep)^\frac12, T-(2M\ep)^\frac12)$. 
In particular, $u^\ep$ and $v_\ep$ have almost everywhere a second order expansion and at
such a point $(x,t)\in \R^N\times ((2M\ep)^\frac12, T-(2M\ep)^\frac12)$ one has
$$
-\tilde u^\ep_t(x,t)-\delta M^-[\tilde u^\ep(\cdot,t)](x) +\delta\mu^{1-q}|D\tilde u^\ep(x,t)|^q +\delta\leq 0
$$
and
$$
-v^\ep_t(x,t)-\delta M^-[v^\ep(\cdot,t)](x) +\delta|Dv^\ep(x,t)|^q +\delta\geq 0\;.
$$
As usual we prove (\ref{muuv}) by contradiction and assume that $\sup_{x,t} \tilde u(x,t)-v(x,t)>0$. 
Then, for any $\alpha>0$, $\sigma>0$ sufficiently small, one can choose $\ep>0$ sufficiently small such that
the map $(x,t)\to \tilde u^\ep(x,t)-v_\ep(x,t)-\alpha|x|^2+\sigma t$ reaches its maximum on $\R^N\times [0,T]$ 
at some point $(\bar x,\bar t)\in \R^N\times ((2M\ep)^\frac12, T-(2M\ep)^\frac12)$. For any $\eta>0$, the point $(\bar x,\bar t)$ is a strict
maximum of the map $(x,t) \to \tilde u^\ep(x,t)-v_\ep(x,t)-\alpha|x|^2+\sigma t-\eta(|x-\bar x|^2+(t-\bar t)^2)$. Jensen's Lemma (see \cite{CIL}) 
then states that one can find $p^n=(p^n_x,p^n_t)\in \R^{N+1}$ such that $p^n\to 0$ and the map 
$$
(x,t) \to \tilde u^\ep(x,t)-v_\ep(x,t)-\alpha|x|^2+\sigma t-\eta (|x-\bar x|^2+(t-\bar t)^2) - \lg p^n_x,x\rg -p^n_tt
$$
has a maximum at some point $(x_n,t_n)$ where $\tilde u^\ep$ and $v_\ep$
have a second order expansion. Note that $(x_n,t_n)\to (\bar x,\bar t)$. At the point $(x_n,t_n)$ we have 
$$
\tilde u^\ep_t(x_n,t_n)=v^\ep_t(x_n,t_n)-\sigma+2\eta(t_n-\bar t)+p^n_t\;,
$$
$$
D\tilde u^\ep(x_n,t_n)=Dv^\ep(x_n,t_n)+2\alpha x_n+2\eta(x_n-\bar x)+p^n_x\;,
$$
\be\label{eq1}
-\tilde u^\ep_t(x_n,t_n)-\delta M^-[\tilde u^\ep(\cdot,t_n)](x_n) +\delta\mu^{1-q}|D\tilde u^\ep(x_n,t_n)|^q +\delta\leq 0\;,
\ee
and
\be\label{eq2}
-v^\ep_t(x_n,t_n)-\delta M^-[v^\ep(\cdot,t_n)](x_n) +\delta|Dv^\ep(x_n,t_n)|^q +\delta\geq 0\;. 
\ee
For the optimality conditions, we have, for any $x\in \R^N$,
$$
\tilde u^\ep(x,t_n)\leq v_\ep(x,t_n)+\tilde u^\ep(x_n,t_n)-v_\ep(x_n,t_n)+\alpha(|x|^2-|x_n|^2)
+\eta (|x-\bar x|^2-|x_n-\bar x|^2) - \lg p^n_x,x-x_n\rg
$$
so that
$$
M^-[\tilde u^\ep(\cdot,t_n)](x_n) \leq M^-[v^\ep(\cdot,t_n)](x_n)+2\alpha +2\eta\;.
$$
Let us set $\xi_n=Dv^\ep(x_n,t_n)$ and estimate the difference between (\ref{eq1}) and (\ref{eq2}): we get
$$
\sigma-2\eta(t_n-\bar t)-p^n_t-2\alpha -2\eta +\delta\left(\mu^{1-q} |\xi_n+2\alpha x_n+2\eta(x_n-\bar x)+p^n_x|^q-|\xi_n|^q\right) \leq 0\;.
$$
When $n\to+\infty$, $\xi_n$ remains bounded since $v^\ep$ is semi-concave. So we can assume that $\xi_n\to \xi_{\alpha,\ep}$ with
$$
\sigma-2\alpha -2\eta +\delta\left(\mu^{1-q} |\xi_{\alpha,\ep}+2\alpha \bar x|^q-|\xi_{\alpha,\ep}|^q\right) \leq 0\;.
$$
Since, $\tilde u$ and $v$ are bounded, so is $\alpha |\bar x|^2$. So $\alpha \bar x$ is bounded (in fact $\alpha \bar x\to 0$ as
$\alpha\to 0$). Then the above inequality implies that $\xi_{\alpha,\ep}$ is bounded, because since $\mu<1$ and $q>1$.  
So, letting $\eta,\ep\to 0$ and then 
$\alpha\to 0$, we get that $\sigma\leq 0$, which contradicts our assumption on $\sigma$. 
\hfill \QED

\begin{Lemma} 
\label{Yn}
Let $u\in {\mathcal C}(\R^N\times [0,T])$ be a supersolution of (\ref{HJnl+})
satisfying $|u|\leq M$ in $\R^N\times (0,T)$. Then, for all $(\bar x,\bar t)\in\R^N\times [0,T]$, $n\in\N$, $R>0$ large and $\sigma>0$ small, there exist a $(\FR_{\bar t,s})$-adapted c\`adl\`ag process $Y^n$ and a control $\zeta^n\in L^p_{\rm ad}([\bar t,T])$ such that
\be\label{Ynxinnlgen}
u(\bar x,\bar t) \geq \E\left[ u(Y^n_{t},t) + C_{_+} \int_{\bar t}^{t}\left|\zeta^n_s \right|^pds-(\delta+\tau) (t-\bar t)-c_n(\sigma,R) \right]\qquad \forall 
t\in [\bar t,T]\; ,
\ee
where
$$
c_n(\sigma,R)= C\left( R^{-p} + \frac{\tau^{p-1}}{\sigma^p}+\omega(\sigma) + (\delta+\tau)\tau\right)\; ,
$$
with $\omega$ the modulus of continuity of $u$ in $B_R(\bar x)\times [0,T]$, and $C$ an universal constant.
\end{Lemma}

\noindent {\it Proof: } 
For any 
$(y,t)\in \R^N\times [0,T)$, $\zeta\in L^p_{\rm ad}( [t,T])$, $a\in {\mathcal A}(t)$, let us denote by
$Y^{x, t, \zeta,a}$ the solution to
$$
\left\{\begin{array}{l}
dY_s= \zeta_sds+ dM^a_s, \; t\leq s\leq T,\\
Y_{t}=x,
\end{array}\right.
$$
where the martingale $M^a$ is defined by (\ref{DefMa}). \\
Let us now fix an initial condition $(\bar x,\bar t)\in\R^N\times [0,T]$. 
For a large $n\in\N$, 
we set
$$
\tau=(T-\bar t)/n\qquad {\rm and }\qquad t_k=\bar t+ k\tau \qquad {\rm for}\; k\in \{0, \dots,n\}\;.
$$
We will build some controls $\zeta^n\in L^p_{\rm ad}( [\bar t,T])$ and $a^n\in {\mathcal A}(\bar t)$ such that
the process $Y^n= Y^{\bar x,\bar t, \zeta^n,a^n}$ satisfies the relation
\begin{equation}\label{Ynxinnl}
u(\bar x,\bar t) \geq \E\left[ u(Y^n_{t_{k}},t_{k}) + C_{_+} \int_{\bar t}^{t_{k}}\left|\zeta^n_s \right|^pds-(\delta+\tau) (t_k-\bar t) \right]\qquad \forall 
k\in\{1,\dots, n\}\;,
\end{equation}
and then deduce from (\ref{Ynxinnl}) that $(Y^n,\zeta^n)$ also satisfy (\ref{Ynxinnlgen}).
We follow closely the construction in \cite{CC}. \\
For any $k\in \{1,\dots, n\}$, let $v^k$ be the solution of (\ref{HJnl+}), defined on the time interval $[0,t_k]$, with terminal condition $u(\cdot, t_k)$. 
 From the representation formula given in Proposition \ref{prop:repHJnl+} we have, for all $x\in\R^N$,
$$
v^k(x,t_{k-1})=\inf_{(\zeta,a)\in L^p_{\rm ad}([t_{k-1},t_k])\times   {\mathcal A}(t_{k-1})}
\E\left[ \ u(Y_{t_k}^{x,t_{k-1},\zeta,a},t_k)+ C_{_+}\int_{t_{k-1}}^{t_k} |\zeta_s|^pds -\delta\tau\right]\; .
$$
Since the filtration $(\FR_{t_{k-1},s})$ is generated by a random Poisson measure, the set $L^p_{\rm ad}([t_{k-1},t_k])\times   {\mathcal A}(t_{k-1})$ is a complete separable space. Moreover $u(\cdot, t_k)$ is continuous. Therefore, thanks to the measurable selection theorem (see \cite{aufr}),   one can build 
Borel measurable maps $x\to Z^{x,k}$ and $x\to A^{x,k}$ from $\R^N$ to $L^p_{\rm ad}( [t_{k-1}, t_k])$
and ${\mathcal A}(t_{k-1})$ respectively, such that
\begin{equation}
\label{vy}
v^k(x,t_{k-1})\geq  \E\left[ u(Y_{t_k}^{x,t_{k-1}, Z^{x,k},A^{x,k}},t_k)+
C_{_+}\int_{t_{k-1}}^{t_k}|Z^{x,k}_s|^pds-(\delta +\tau)\tau\right]\qquad \forall x\in \R^N\;.
\end{equation}
We now construct  $\zeta^n$, $a^n$ and  $Y^n$ by induction on the time intervals $[t_{k-1}, t_{k})$:\\
 On $[\bar t, t_1)$ we set
$\zeta^n_t=Z^{\bar x,1}_t$, $a^n=A^{\bar t,1}$ and  $Y^n= Y^{\bar x,\bar t,\zeta^n,a^n}$. Assume that $\zeta^n$, $Y^n$ and $a^n$
have been built on $[\bar t, t_{k-1})$. Then we set 
$$
\zeta^n= Z^{Y^n_{t_{k-1}},k}, \; a^n= A^{Y^n_{t_{k-1}},k}, 
 \qquad {\rm and }\qquad Y^n= Y^{\bar x,\bar t,\zeta^n,a^n}\qquad  {\rm on }\; [t_{k-1}, t_k)\;.
$$
(The process $Y^{\bar x,\bar t,\zeta^n,a^n}$ is $P$-a.s. continuous on each fixed $t$, so we have $Y^n_{t_{k-1}-}=Y^n_{t_{k-1}}$ $P$-a.s., which means that  $\zeta^n$ and $A^n$ are defined $P$-as surely. )\\ 
We remark that, on $[t_{k-1}, t_k)$, we have 
$Y^n=Y^{Y^n_{t_{k-1}},t_{k-1},\zeta^n,a^n}$.\\
Let us fix now some $k\in\{ 1,\ldots,n\}$.
Since the processes $A^{x,k}$ and $Z^{x,k}$ are $(\FR_{t_{k-1},s})$-adapted and therefore independent of $\FR_{\bar t,t_{k-1}}$, the same holds also for $M^{A^{x,k}}_{t_k}$ and finally for $Y_{t_k}^{x,t_{k-1}, Z^{x,k},A^{x,k}}$, while $Y^n_{t_{k-1}}$ is $\FR_{\bar t,t_{k-1}}$-measurable. It follows that
\[ \E\left[ u(Y_{t_k}^{x,t_{k-1}, Z^{x,k},A^{x,k}},t_k)+
C_{_+}\int_{t_{k-1}}^{t_k}|Z^{x,k}_s|^pds\right]_{x=Y^n_{t_{k-1}}}
=\E\left[ u(Y_{t_k}^n,t_k)+
C_{_+}\int_{t_{k-1}}^{t_k}|\zeta^n_s|^pds\;\Big\vert\;\FR_{\bar t,t_{k-1}}\right].\]
Using (\ref{vy}), the fact that $u$ is a supersolution of (\ref{HJnl+}) and the comparison Lemma \ref{compar}, this leads to the relation
$$
u(Y^n_{t_{k-1}},t_{k-1}) \geq \E\left[ u(Y^n_{t_{k}},t_{k}) + C_{_+} \int_{t_{k-1}}^{t_{k}}\left|\zeta^n_s \right|^pds-(\delta+\tau) \tau\; \big|\; {\cal F}_{\bar t,t_{k-1}}\right]
\qquad \P-{\rm a.s.}\;.
$$
Taking the expectation on both sides of the above inequality and summing up gives \eqref{Ynxinnl}.\\

\noindent We now extend this inequality to the full interval $[\bar t, T]$ and prove (\ref{Ynxinnlgen}).
Let $t\in [\bar t,T]$ and $k$ be such that $t\in [t_{k-1},t_k)$. From \eqref{Ynxinnl}, we have
\be\label{preuve:Ynxinnlgen}
\begin{array}{r}
u(\bar x,\bar t) \geq \E\left[ u(Y^n_{t},t) + C_{_+} \int_{\bar t}^{t}\left|\zeta^n_s \right|^pds-(\delta+\tau) (t-\bar t)\right]\qquad\qquad\qquad\qquad\\ 
+\E\left[ u(Y^n_{t_{k-1}})-u(Y^n_t,t)\right] -(\delta+\tau)\tau 
\end{array}
\ee
Let us fix $R>0$ and $\sigma>0$.
Since $u$ is bounded by $M$ and from the definition of the modulus $\omega$, we have
\begin{equation}
\begin{array}{rl}
\label{preuve:pp}
\E\left[ u(Y^n_{t_{k-1}})-u(Y^n_t,t)\right]\leq &
\omega(\sigma)\P\left[|Y^n_{t_k}-\bar x|\leq R, \ |Y^n_{t}-\bar x|\leq R,\ |Y^n_{t_k}-Y^n_{t}|\leq \sigma\right]\\
&+2 M(\P\left[|Y^n_{t_k}-\bar x|> R\right]+\P\left[|Y^n_{t}-\bar x|> R\right]+
\P\left[ |Y^n_{t_k}-Y^n_{t}|> \sigma\right]) 
\end{array}
\end{equation}
To estimate the right hand side term of (\ref{preuve:pp}), we first note  that, for any $0<s<t$, it holds that
$$
\E\left[ \left|Y^n_t-Y^n_s\right|^p\right]
\leq 2^{p-1} \left\{ \E\left[ \left|\int_s^t \zeta^n_\tau d\tau \right|^p\right]
+\E\left[ \left|M^{a_n}_t-M^{a_n}_s\right|^p\right] \right\}
$$
But, thanks to (\ref{Ynxinnl}) again, we have
\be\label{BorneZetanl}
\E\left[\int_{\bar t}^T \left|\zeta^n_\tau\right|^p d\tau \right]\; \leq \; 2M+(\delta+\tau)T\; \leq\;  C
\ee
so that, by H\"{o}lder's inequality, 
$$
\E\left[ \left|\int_s^t \zeta^n_\tau d\tau \right|^p\right] \leq C(t-s)^{p-1}\;.
$$
Also by H\"{o}lder we have 
$
\E\left[ \left|M^{a_n}_t-M^{a_n}_s\right|^p\right]\leq \left( \E\left[ \left|M^{a_n}_t-M^{a_n}_s\right|^2\right]\right)^{p/2}
$
where, by It\^{o},  
\be\label{estiMn}
\begin{array}{rl}
\ds{  \E\left[ \left|M^{a_n}_t-M^{a_n}_s\right|^2\right] \; = } & {\ds  \E\left[\int_s^t \int_{\BN} \lambda_s^2|b_s|^2{\bf 1}_{A_s}(e)d\nu(e)ds \right] }\\
 & {\ds  =\delta \E\left[ \int_s^t \lambda_s^2 ds \right]\; \leq\; \delta(t-s)}
\end{array}
\ee
To summarize
$$
\E\left[ \left|Y^n_t-Y^n_s\right|^p\right]
\; \leq \; C ((t-s)^{p-1}+(t-s)^{p/2}) \; \leq \; C(t-s)^{p-1}
$$
since $p<2$. Therefore we get
$$
\P[ |Y^n_{t_{k}}-\bar x|>R]+\P[ |Y^n_{t}-\bar x|>R]+ \P[ |Y^n_t-Y^n_{t_{k}}|>\sigma]
\leq C ( R^{-p} + \frac{|t-t_k|^{p-1}}{\sigma^p})
$$
which, coming back to (\ref{preuve:Ynxinnlgen}) and (\ref{preuve:pp}), proves claim (\ref{Ynxinnlgen}). \hfill \QED

\begin{Lemma}\label{lem:IneqStochnl} Let $u\in {\mathcal C}(\R^N\times [0,T])$ be a supersolution of (\ref{HJnl+})
satisfying $|u|\leq M$ in $\R^N\times (0,T)$. Then, for any $(\bar x,\bar t)\in \R^N\times (0,T)$ there is a stochastic
basis $(\bar\Omega,\bar{\cal F}, \bar\P)$, a filtration $(\bar{\cal F}_t)_{t\geq \bar t}$, a c\`{a}dl\`{a}g process $(\bar Y_t)$ adapted to $(\bar{\cal F}_t)_{t\geq \bar t}$
and a process $\bar\zeta\in L^p_{\rm ad}( [\bar t,T])$ such that 
\begin{equation}\label{IneqStochnl}
u(\bar x,\bar t)\geq \E\left[ u(\bar Y_{t},t) +C_{_+}\int_{\bar t}^t|\bar\zeta_s|^pds\right] -\delta (t-\bar t)\qquad \forall t\in (\bar t, T)\;,
\end{equation}
where $C_{_+}>0$ is the  universal constant given by 
(\ref{eq:C+}), 
and 
\be\label{ineq:stochprocnl}
\E\left[ \left| \bar Y_t-\bar x-\int_{\bar t}^t \bar\zeta_sds \right|^r\right] \leq \delta^{\frac r2}|t-\bar t|^{r/2}\qquad \forall t\in [\bar t, T)
\ee
for any $r\in(0,2]$. 
\end{Lemma}

\noindent {\it Proof: } 
This  Lemma will follow from Lemma \ref{Yn} by passing to the limit as $n\to+\infty$ in (\ref{Ynxinnlgen}). For this we set
$$
\Lambda^n_t\dot{=}\int_{\bar t}^t \zeta^n_sds \qquad \forall t\in[\bar t, T]\;.
$$
From (\ref{BorneZetanl}), the sequence of probability measures
$(\P_{\Lambda^n})$ on ${\mathcal C}([\bar t,T], \R^{N})$ is tight. Let ${\bf D}(\bar t)$ be the set of c\`{a}dl\`{a}g
functions from $[\bar t,T]$ to $\R^N$, endowed with the Meyer-Zheng topology (see \cite{mezh}). Since $\E[|M^{a_n}_T|]$ is uniformly bounded
(thanks to (\ref{estiMn})), Theorem 4 of \cite{mezh} states that the sequence of martingale measures $(\P_{M^{a_n}})$ is
tight on ${\bf D}(\bar t)$. Then, from Prohorov's Theorem (Theorem 4.7 of \cite{KS}), we can find a subsequence of $(Y^n,\Lambda^n)$,
again labeled $(Y^n,\Lambda^n)$, and a measure $m$ on ${\mathcal C}([\bar t,T], \R^{N})\times {\bf D}(\bar t)$ such that 
$(\P_{(Y^n,\ \Lambda^n)})$ weakly converges to $m$. Skorokhod's embedding Theorem (Theorem 2.4 of \cite{DPZ}) implies that we can find 
random variables $(\bar Y^n,\bar \Lambda^n)$ and $(\bar Y,\bar \Lambda)$ 
defined on a new probability space $(\bar \Omega, \bar {\cal A}, {\bf \bar P})$, such that
$(\bar Y^n,\bar \Lambda^n)$ has the same law as $(Y^n,\Lambda^n)$ for any $n$,  the law of $(\bar Y,\bar \Lambda)$ is $m$
and,  $\bar \P-$almost surely,
the sequence $(\bar \Lambda^n)$ converges to $(\bar \Lambda)$ in
${\mathcal C}([\bar t,T], \R^{N})$ while, for any $t$ belonging to some set $I\subset [\bar t,T]$ 
of full measure in $[\bar t,T]$, the sequence $(\bar Y^n_t)$ converges to $\bar Y_t$ (Theorem 5 of \cite{mezh}). 

Since $t\to \Lambda^n_t$ is absolutely continuous $\P-$a.s. and since $\bar \Lambda^n$ has the same law as $\Lambda^n$,
$t\to \bar \Lambda^n_t$ is absolutely continuous $\bar \P-$a.s.. Let us set $\bar \zeta^n_s= \frac{d}{ds} \bar \Lambda^n_s$. Then, 
by \eqref{BorneZetanl}, 
$
\E[\int_{\bar t}^T \left|\bar \zeta^n_s \right|^pds]\leq C
$ for all $n\geq 0$.
Therefore,  up to a subsequence again labeled in the same way, $(\bar \zeta^n)$ converges weakly in $L^p( [\bar t, T])$ to some 
limit, $\bar \zeta$, which, $\bar \P-$a.s., satisfies $\bar \Lambda_t=\int_{\bar t}^t\bar \zeta_sds$
for all $t\in [\bar t, T]$. 

Note that $\bar M^n_t\dot{=} \bar Y^n_t-\bar x-\bar \Lambda^n_t$
has the same law as $M^{a_n}_t$, so that by H\"{o}lder and (\ref{estiMn}), 
$$
\bar \E\left[ \left|\bar M^n_t\right|^r\right] \leq   \delta^{\frac r2}(t-\bar t)^{r/2}
$$
for all $r\in(0,2]$ and for all $t\in [\bar t,T]$. Passing to the limit in the above inequality gives
$$
\bar \E\left[ \left|\bar Y_t-\bar x-\bar \Lambda_t\right|^r\right]\leq\delta^{\frac r2} (t-\bar t)^{r/2} \qquad \forall t\in I,\; \forall r\in(0,2]\;.
$$
We get the above inequality for all $t\in [\bar t, T)$ thanks to the c\`{a}d\`{a}g property of the trajectories
of $Y$. Recalling \eqref{Ynxinnlgen}, a classical lower semicontinuity argument yields
$$
u(\bar x,\bar t) \geq \bar \E\left[ u(\bar Y_{t},t) + C_{_+} \int_{\bar t}^{t}\left|\bar \zeta_s \right|^pds-\delta (t-\bar t) \right]\qquad \forall t\in I\;,
$$
and we conclude the proof by using again the c\`{a}d\`{a}g property of the trajectories
of $Y$.  \hfill \QED

\section{Analysis of subsolutions}

In this section we investigate properties of subsolution of equation (\ref{HJnl-}). 

\begin{Proposition}\label{prop:SuperSol-nl} For any fixed $(x,t)\in \R^N\times (0,T]$, there is a continuous supersolution 
 $w$ of (\ref{HJnl-}) in $\R^N\times [0,t)$ such that
 \be\label{estiw}
 \frac{1}{C} (t-s)^{1-p}|x-y|^p - C (t-s)^{1-p/2} \leq w(y,s) \leq  C (t-s)^{1-p}|x-y|^p + C (t-s)^{1-p/2}
 \ee
for any $(y,s)\in \R^N\times [0,t)$ and for some universal constant $C$.
 \end{Proposition}

{\bf Proof : } It relies on control interpretation of equation (\ref{HJnl-}) as well as the construction of some Brownian bridges (see \cite{FY}). 
Let us assume, without loss of generality, that $x=0$. 
Having fixed $\alpha\in (1-1/p,1/2)$, $(y,s)\in \R^N\times [0,t)$ and $a\in {\mathcal A}(s)$, let 
$Y_\tau^{y,s,a}$ be the solution to 
$$
\left\{\begin{array}{l}
dY_\tau= -\alpha\,\frac{Y_\tau}{t-\tau}\,d\tau+  dM^a_\tau\\
Y_s=y
\end{array}\right.
$$
Then one easily checks that
\begin{equation}\label{FormuleY}
Y_\tau= (t-s)^{-\alpha}(t-\tau)^{\alpha}y+(t-\tau)^{\alpha}\int_s^\tau (t-\sigma)^{-\alpha} dM^a_\sigma\;. 
\end{equation}
Let us set 
$$
Z^{y,s,a}_\tau \doteq  -\,\alpha\,Y_\tau/(t-\tau) \qquad {\rm and }\qquad J(y,s,a)=\E\left[\ \int_s^t \left|Z^{y,s,a}_\tau\right|^pd\tau \ \right]\;.
$$
We claim that there is a universal constant $C>0$ such that 
\be\label{estiJ}
 \frac{1}{C} (t-s)^{1-p}|x-y|^p - C (t-s)^{1-p/2} \leq J(y,s,a) \leq  C (t-s)^{1-p}|x-y|^p + C (t-s)^{1-p/2}\;.
 \ee
 Indeed
\begin{eqnarray*}
\lefteqn{J(y,s,a)\ =\ \E\left[\int_s^t |Z^{y,s,a}_\tau|^pd\tau\right]} \\
\; & \leq & 2^{p-1} \alpha ^p(t-s)^{-\alpha p}|y|^p \int_s^t (t-\tau)^{p(\alpha-1)}d\tau 
\\
& & \hspace{3cm}+2^{p-1} \alpha^p\int_s^t (t-\tau)^{p(\alpha-1)}\E\left[\Big| \int_s^\tau (t-\sigma)^{-\alpha}dM^a_\sigma \Big|^p\right]d\tau
\end{eqnarray*}
where, by H\"{o}lder and It\^{o}, 
$$
\E\left[\Big| \int_s^\tau (t-\sigma)^{-\alpha}dM^a_\sigma \Big|^p\right]
\leq C(t-s)^{\frac p2(1-2\alpha)}  
$$
So 
$$
J(y,s,a)\leq C (t-s)^{1-p}|y|^p + C(t-s)^{1-p/2}\;.
$$
In the same way, 
\begin{eqnarray*}
\lefteqn{J(y,s,a)\ =\ \E\left[\int_s^t |Z^{y,s,a}_\tau|^pd\tau\right]} \\
\; & \geq & 2^{1-p} \alpha ^p(t-s)^{-\alpha p}|y|^p \int_s^t (t-\tau)^{p(\alpha-1)}d\tau 
\\
& & \hspace{3cm}-\alpha^p\int_s^t (t-\tau)^{p(\alpha-1)}\E\left[\Big| \int_s^\tau (t-\sigma)^{-\alpha}dM^a_\sigma\Big|^p\right]d\tau\\
&\geq & (1/C) (t-s)^{1-p}|y|^p - C(t-s)^{1-p/2}\,,
\end{eqnarray*}
Whence (\ref{estiJ}). 

Next we introduce the value function $w$ of the optimal control problem
$$
w(y,s)=  C_{_-}\ \sup_{a\in {\mathcal A}(s)}  J(y,s,a)- \delta (t-s)\;,
$$
where $C_{_-}\doteq\frac{\delta^{p/q}}{pq^{p/q}}$.
Let us first show that $w$ is continuous on $\R^N\times [0,t)$. 
The map $y\to J(y,s,a)$ being convex (since the map $y\to Y^{y,s,a}_\tau$ is affine and $p>1$) and locally uniformly bounded
(thanks to (\ref{estiJ})), it is has a modulus of continuity which is locally
uniform with respect to $s$ and $a$. The map $s\to J(y,s,a)$
being locally H\"{o}lder continuous on $[0,t)$, locally uniformly with respect to $y$ and $a$, this implies that the map
$(y,s) \to J(y,s,a)$ has a modulus of continuity which is uniform with respect to $a$. Therefore $w$ is continuous
on $\R^N\times [0,t)$. 

Using the fact that $w$ is continuous and arguments similar to the ones in \cite{BMR} one can prove that $w$ satisfies the Hamilton-Jacobi equation 
$$
\begin{array}{l}
\ds{  -w_t+ \inf_{\lambda\in(0,1],\ b\in \BN\backslash \{0\}} \left\{- \lg  -\alpha\,\frac{y}{t-s}, Dw\rg + C_{_-}\left|-\alpha\,\frac{y}{t-s}\right|^p \right.  }\\
\ds{  \qquad \qquad \qquad \qquad \left.  -\delta\  \frac{w(y+\lambda b,s)-w(y,s)-\lg Dw(y,s),\lambda b\rg}{|b|^2} \right\}-\delta=0  }
 \end{array}
$$
Since
$$
 \inf_{\lambda\in(0,1],\ b\in \BN\backslash \{0\}}\left\{- \ \frac{w(y+\lambda b,s)-w(y,s)-\lg Dw(y,s),\lambda b\rg}{|b|^2}\right\} = - M^+[w(\cdot,s)](y)
 $$
 while 
$$
- \lg -\alpha\,\frac{y}{t-s}, Dw\rg +  C_{_-}\left|-\alpha\,\frac{y}{t-s}\right|^p \geq \frac{1}{\delta} |Dw|^q\;,
$$
$w$ is a supersolution of (\ref{HJnl-}). We finally note that $w$ satisfies (\ref{estiw}) because the inequalities 
(\ref{estiJ}) are uniform with respect to $a$. 
\hfill\QED

\begin{Lemma}\label{lem:EstiSubSolStoch} Let $u\in {\mathcal C}(\R^N\times [0,T])$ be a subsolution of (\ref{HJnl-})
satisfying $|u|\leq M$.
Then, 
for all $(x,t)\in \R^N\times (0,T)$  and all $(y,s)\in \R^N \times [0,t)$,  
\begin{equation}\label{eq:EstiSubSolStoch}
u(y,s)  \leq  u(x,t)+C\left\{ |y-x|^p(t-s)^{1-p}+ (t-s)^{1-p/2}\right\}
\end{equation}
for some universal constant $C>0$.
\end{Lemma}

\begin{Remark}{\rm In particular, if $u=u(x)$ is a subsolution of the stationary equation
$$
-\delta M^+[u(\cdot,t](x) +\frac{1}{\delta}|Du|^q -\delta=0 \qquad {\rm in }\quad\R^N,
$$
then inequality (\ref{eq:EstiSubSolStoch}) implies that, for any $x,y\in \R^N$ and any  $\tau>0$, 
$$
u(x)  \leq  u(y)+C\left\{ |y-x|^p\tau^{1-p}+\tau^{1-p/2}\right\}\,,
$$
for some universal constant $C$. Thus, choosing $\tau=|x-y|^2$ yields
$
u(x)  \leq  u(y)+C\ |y-x|^{2-p}
$, that is, $u$ is H\"{o}lder continuous. This extends to nonlocal equations one of the results of \cite{CDLP}. 
}\end{Remark}

\noindent {\it Proof:} According to Proposition \ref{prop:SuperSol-nl} there is a supersolution  $w$ of (\ref{HJnl-})  which satisfies
$$ 
\frac{1}{C} (t-s)^{1-p}|x-y|^p - C (t-s)^{1-p/2} \leq w(y,s) \leq  C (t-s)^{1-p}|x-y|^p + C (t-s)^{1-p/2}
$$
for any $(y,s)\in \R^N\times [0,t)$ and for some universal constant $C$. Since $u$ is continuous and bounded and $u(y,t) \leq \lim_{s\to t} w(y, s)+u(x,t)$
for any $y\in \R^N$,  we get $u\leq w+ u(x,t)$  on $\R^N\times [0,t)$ by comparison (Lemma \ref{compar}). Whence the result. 
\hfill \QED

\begin{Lemma}\label{lem:EstiSubSolStoch2} Let $u\in {\mathcal C}(\R^N\times [0,T])$ be a subsolution of  \eqref{HJnl-} satisfying $|u|\leq M$. Fix
$(\bar x,\bar t)\in \R^N\times (0,T)$, $\zeta\in L^p_{\rm ad}([\bar t, T])$ and let $(X_t, \zeta_t)$ be stochastic processes satisfying
(\ref{ineq:stochprocnl}). Then, for any $x\in \R^N$ and $t\in (\bar t, T)$,
\begin{equation}\label{eq:EstiSubSolStoch2}
\begin{array}{l}
u(x,\bar t)  - \E[u(X_t,t)]\\
\qquad \leq C\left\{(t-\bar t)^{1-p}\left(\E\left[\left(\int_{\bar t}^t |\zeta_s|ds\right)^p\right]+ 
 |\bar x-x|^p\right)+ (t-\bar t)^{1-p/2} \right\}
\end{array}
\end{equation}
for some constant $C>0$. 
\end{Lemma}
{\it Proof:} 
Fix $t\in (\bar t,T)$ and apply Lemma \ref{lem:EstiSubSolStoch} to $(x,\bar t)$
and $(X_{t}(\omega),t)$. Then, for almost all $\omega\in \Omega$,
$$
u(x,\bar t) \; \leq \; u(X_t(\omega),t)+C\left\{ |X_t(\omega)-x|^p(t-\bar t)^{1-p}+ (t-\bar t)^{1-p/2}\right\}\;.
$$
Hence,
$$
u(x,\bar t)\; 
\leq \; \E\left[ u(X_t,t)\right]+
C\left\{ (\E\left[|X_t-\bar x|^p\right]+|\bar x-x|^p)(t-\bar t)^{1-p}+ (t-\bar t)^{1-p/2}\right\}\;.
$$
Since, on account of (\ref{ineq:stochprocnl}),
$$
\E\left[|X_t-\bar x|^p\right]\leq C\left\{ \E\left[\left(\int_{\bar t}^t |\zeta_s|ds\right)^p\right]+
 (t-\bar t)^{\frac{p}{2}} \right\},
$$
the conclusion follows. \hfill \QED

In order to proceed, we need to recall the following weak reverse H\"{o}lder inequality:

\begin{Lemma}[\cite{CC}]\label{lem:RevHol}
Let $(\Omega, {\cal A}, \P)$ be a probability space. Let $p\in (1,2)$ and assume that the function
$\xi\in L^p( (a,b);\R^N)$ satisfies the inequality
\begin{equation}\label{RevHolStoch}
\E\left[ \frac{1}{t-a} \int_{a}^t |\xi_s|^pds\right] \leq 
A\,\E\left[\Big(\frac{1}{t-a}\int_{a}^t|\xi_s|ds\Big)^p\right]+\frac{B}{(t-a)^{\frac p2}}\qquad \forall t\in (a,b]
\end{equation}
for some positive constants $A$ and $B$. 
Then there are constants $\theta\in(p,2)$ and $C> 0$, depending only on $p$ and $A$, such that
$$
\E \left[\Big(\int_{a}^t|\xi_s|ds\Big)^p\right] \leq C (t-a)^{p-\frac{p}{\theta}}\left\{(b-a)^{\frac{p}{\theta}-1}\|\xi\|^p_p 
+B(b-a)^{\frac{p}{\theta}-{\frac p2}}\right\}
\qquad \forall t\in (a,b]\;.
$$
\end{Lemma}

Thanks to this inequality, we can estimate the $L^p$ norm of the process $\zeta$ appearing in Lemma \ref{lem:IneqStochnl}. 

\begin{Lemma}\label{LemABStoch} Let $u\in {\mathcal C}(\R^N\times [0,T])$ be a subsolution of \eqref{HJnl-}
such that $|u|\leq M$ and let $\tau\in (0,T)$.
Then there is a universal constant $\theta\in(p,2)$ and a constant $C(\tau,\delta)>0$ such that, for every 
$(\bar x,\bar t)\in \R^N\times (0,T-\tau)$, 
and every stochastic processes $(X_t, \zeta_t)$ satisfying (\ref{IneqStochnl}) and (\ref{ineq:stochprocnl}), we have
$$
\E \left[\left(\int_{\bar t}^t|\zeta_s|ds\right)^p\right] 
\leq C(\tau)(t-\bar t)^{p-\frac{p}{\theta}}\qquad \forall t\in (\bar t, T)\;.
$$
\end{Lemma}
{\it Proof:} First, observe that, by Lemma \ref{lem:EstiSubSolStoch2} applied to $x=\bar x$, 
$$
u(\bar x,\bar t)  \leq  \E[u(X_{t},t)]
+ C\left((t-\bar t)^{1-p}\E\left[\left(\int_{\bar t}^t |\zeta_s|ds\right)^p\right]+  (t-\bar t)^{1-p/2} \right)
$$
for all $t\in[\bar t,T)$. Moreover, in view of (\ref{IneqStochnl}),
$$
\E\left[ u(X_{t},t)\right] \leq u(\bar x,\bar t)-  C_{_+} \E\left[ \int_{\bar t}^t|\zeta_s|^pds\right]+\delta(t-\bar t) \qquad \forall t\in [\bar t, T)\;.
$$
Hence, taking into account that $t-\bar t\leq C(t-\bar t)^{1-p/2}$,
$$
\E\left[ \int_{\bar t}^t|\zeta_s|^pds\right]  
\leq C\, (t-\bar t)^{1-p}\E\left[\left(\int_{\bar t}^t |\zeta_s|ds\right)^p\right]+ C(t-\bar t)^{1-p/2}\qquad \forall t\in [\bar t, T)\;.
$$
Then, owing to Lemma \ref{lem:RevHol}, there
are universal constants $\theta\in(p,2)$ and $C>0$ such that 
$$
\E \left[\left(\int_{\bar t}^t|\zeta_s|ds\right)^p\right] \leq C\left(\|\zeta\|_p^p+1\right)
\frac{ (t-\bar t)^{p-\frac{p}{\theta}}}{(T-\bar t)^{1-\frac{p}{\theta}}}
\qquad \forall t\in (\bar t,T)\;.
$$
Since $u$ is bounded by $M$, assumption (\ref{IneqStochnl}) implies that $\|\zeta\|_p\leq C$. 
So, we finally get 
$$
\E \left[\left(\int_{\bar t}^t|\zeta_s|ds\right)^p\right] \leq C(\tau)(t-\bar t)^{p-\frac{p}{\theta}}\qquad \forall t\in (\bar t, T)\;,
$$
because $\bar t\leq T-\tau$.
\hfill \QED

\section{Proof of Theorem \ref{Regu2NL}}

Let $u:\R^N\times [0,T]\to \R$ be a continuous  supersolution of  (\ref{HJnl+}), a subsolution of (\ref{HJnl-}) and  such that
$|u|\leq M$. \\

\noindent {\it Space regularity:} Fix $(\bar x,\bar t)\in \R^N\times (0,T-\tau)$ and let $x\in \R^N$.  By Lemma \ref{lem:IneqStochnl} 
there is a control $\zeta\in L^p_{\rm ad}( [\bar t,T])$ and 
an adapted process $X$ such that (\ref{IneqStochnl}) and (\ref{ineq:stochprocnl}) hold. So,
\begin{equation}\label{IneqStochSimple}
u(\bar x,\bar t)\geq \E\left[ u(X_t,t) \right] -\delta (t-\bar t)\qquad \forall t\in [\bar t, T)\;.
\end{equation}
Also, Lemma  \ref{LemABStoch} ensures that 
\begin{equation}\label{EstiZeta4}
\E \left[\left(\int_{\bar t}^t|\zeta_s|ds\right)^p\right] \leq C(\tau)(t-\bar t)^{p-\frac{p}{\theta}}\qquad \forall t\in (\bar t, T)
\end{equation}
for some universal constant $\theta\in (p,2)$ and some constant $C(\tau)>0$. 
Furthermore,  applying Lemma~\ref{lem:EstiSubSolStoch2}, for any $t\in (\bar t, T)$ we have
$$
\begin{array}{l}
u(x,\bar t)  - \E[u(X_{t},t)]\\
\qquad \leq C\left\{ (t-\bar t)^{1-p}\E\left[\left(\int_{\bar t}^t |\zeta_s|ds\right)^p\right]+ 
|\bar x-x|^p(t-\bar t)^{1-p}+ (t-\bar t)^{1-p/2}\right\}\;.
\end{array}
$$
Plugging (\ref{IneqStochSimple}) and (\ref{EstiZeta4}) into the above inequality leads to
$$
u(x,\bar t) \leq u(\bar x,\bar t) +\delta(t-\bar t)+C(\tau)(t-\bar t)^{(\theta-p)/\theta}+C|\bar x-x|^p(t-\bar t)^{1-p}+C(t-\bar t)^{1-p/2}
$$
for any $t\in (\bar t, T)$. \\
Since $1>1-p/2>(\theta-p)/\theta$ (recall that $\theta<2$),
$$
u(x,\bar t) \leq u(\bar x,\bar t) +C(\tau)(t-\bar t)^{(\theta-p)/\theta}+C|\bar x-x|^p(t-\bar t)^{1-p}\;.
$$
Then, for $|x-\bar x|$ sufficiently small,  choose $t=\bar t+|x-\bar x|^{\theta/(\theta-1)}$ to obtain
$$
u(x,\bar t)\; \leq u(\bar x, \bar t)+C(\tau) |x-\bar x|^{(\theta-p)/(\theta-1)}.
$$

\noindent {\it Time regularity : } Let now $t\in (0,T-\tau)$. Then, in light of (\ref{IneqStochSimple}),
$$
 u(\bar x,\bar t) \geq \E\left[ u(X_{t},t) \right]-\delta (t-\bar t).
$$
Now, applying  the space regularity result we have just proved, we obtain
$$
\E\left[ u(X_{t},t) \right]\geq  u(\bar x,t)- C(\tau) \E\left[ \left|X_t-\bar x\right|^{\theta-p\over\theta-1}\right].
$$
Moreover, since $(\theta-p)/(\theta-1)\ <\ 1$, by (\ref{ineq:stochprocnl}) we get
$$
\E\left[|X_t-\bar x|^{\theta-p\over\theta-1}\right] \leq  C \E\left[\Big(\int_{\bar t}^t |\zeta_s|ds\Big)^{\theta-p\over\theta-1}\right]+ 
 C (t-\bar t)^{\theta-p\over 2(\theta-1)}\;.
$$ 
Also, by H\"{o}lder's inequality and (\ref{EstiZeta4}), 
$$
\E\left[\Big(\int_{\bar t}^t |\zeta_s|ds\Big)^{\theta-p\over\theta-1}\right]
\leq C\ \left\{\E\left[\Big(\int_{\bar t}^t |\zeta_s|ds\Big)^p\right]\right\}^{\frac{\theta-p}{p(\theta-1)}}
\leq 
C(\tau) (t-\bar t)^{\theta-p\over\theta}\;.
$$
Notice that $(\theta-p)/(2(\theta-1))> (\theta-p)/\theta$ since $\theta<2$. So,
$$
u(\bar x,\bar t)\geq u(\bar x,t) - C(\tau)  (t-\bar t)^{\theta-p\over\theta}\;.
$$
To derive the reverse inequality, one just needs to apply  Lemma~\ref{lem:EstiSubSolStoch} with $y=x=\bar x$ to get
$$
u(\bar x,\bar t)  \leq  u(\bar x,t) +C(t-\bar t)^{1-p/2}\;.
$$
This leads to the desired result since $1-p/2>(\theta-p)/\theta$. 
\hfill \QED


\begin{thebibliography}{9}
\bibitem{AT}
Alvarez O. and Tourin A.,
{\it  Viscosity solutions of nonlinear integro-differential equations,}
 Ann. Inst. H. Poincaré Anal. Non Lin\'eaire 13(3) (1996) 293--317.

\bibitem{Ar}
Arisawa, M.
{\it  A remark on the definitions of viscosity solutions for the integro-differential equations with Lévy operators.}
  J. Math. Pures Appl. (9)  89  (2008),  no. 6, 567--574. 


\bibitem{aufr}Aubin J.-P.,  Frankowska H.,
                      {\sc Set-Valued Analysis},
                      Birkh\"auser, Boston, 1990.

\bibitem{BBP} Barles G., Buckdahn R.,  Pardoux E., 
{\it Backward stochastic differential equations and integral-partial differential equations.}
  Stochastics Stochastics Rep. 60(1--2) (1997), 57--83.

\bibitem{BI} Barles G., Imbert, C. 
{\it Second-order elliptic integro-differential equations: viscosity solutions' theory revisited.}
  Ann. Inst. H. Poincaré Anal. Non Linéaire  25  (2008),  no. 3, 567--585.

\bibitem{BMR} Buckdahn R., Ma J., and Rainer C. 
{\it Stochastic control problems for systems driven by normal martingales}, 
Ann. Appl. Probab. Volume 18, Number 2 (2008), 632-663. 

\bibitem{CS}
Caffarelli L.,  Silvestre L.
{\it  Regularity theory for fully nonlinear integro-differential equations.}
  Comm. Pure Appl. Math.  62  (2009),  no. 5, 597--638. 

\bibitem{CDLP} Capuzzo Docetta I., Leoni F., Porretta A.,
{\it H\"{o}lder estimates for degenerate elliptic equations with coercive Hamiltonians.} To appear in Transactions Amer. Math. Soc. 

\bibitem{CC} Cannarsa P., Cardaliaguet P.,
{\it H\"{o}lder estimates in space-time for viscosity solutions of Hamilton-Jacobi equations}. To appear in CPAM. 

\bibitem{CIL} Crandall M.G., Ishii H., Lions P.-L.  (1992)
{\it User's guide to viscosity solutions of second order Partial Differential Equations}. 
Bull. Amer. Soc., 27, pp. 1-67.



\bibitem{dll} Da Lio F., Ley O., {\it Uniqueness results for convex Hamilton-Jacobi equations under
$p > 1$ growth conditions on data}, Pre-print. 

\bibitem{DPZ} Da Prato G.,  Zabczyk J., {\sc Stochastic equations in infinite dimensions.}
 Encyclopedia of Mathematics and its Applications, 44. Cambridge University Press, Cambridge, 1992. 


\bibitem{FY} Fleming W. H., Sheu S. J.,
{\it  Stochastic variational formula for fundamental solutions of parabolic PDE.}
 Appl. Math. Optim. 13 (1985), no. 3, 193--204.


\bibitem{KS} Karatzas I., Shreve S. E.,
{\sc Brownian motion and stochastic calculus.} 
Second edition. Graduate Texts in Mathematics, 113. Springer-Verlag, New York, 1991.


 
\bibitem{mezh} Meyer P.A., Zheng W.A., {\it Tightness criteria for laws of semimartingales},
Ann. I.H.P., section B, tome 20, nr.4 (1984), 353-372.

\bibitem{Ph} Pham H., {\it Optimal stopping of controlled jump diffusion processes: A viscosity solution approach.}
 J. Math. Systems Estim. Control. 8(1), 1998. 

\bibitem{Sa1} Sayah A., {\it Equations d'Hamilton-Jacobi du premier ordre avec termes int\'egro-diff\'erentiels. I. Unicit\'e des solutions de viscosit\'e.}
 Comm. Partial Differential Equations, 16(6--7) (1991), 1057--1074. MR1116853 (93b:35027a)

\bibitem{Sa2} Sayah A., {\it Equations d'Hamilton-Jacobi du premier ordre avec termes int\'egro-diff\'erentiels. II. Existence de solutions de viscosit\'e.}
 Comm. Partial Differential Equations, 16(6--7) (1991), 1075--1093.

\bibitem{Si}
Silvestre L.
{\it  H\"{o}lder estimates for solutions of integro-differential equations like the fractional Laplace.}
  Indiana Univ. Math. J.  55  (2006),  no. 3, 1155--1174.

\bibitem{So} Soner H.M., {\it Optimal control of jump-Markov processes and viscosity solutions.}
 Stochastic differential systems, stochastic control theory and applications (Minneapolis, Minn., 1986), 
 IMA Vol. Math. Appl. 10, 501--511, Springer, New York, 1988.

\end{thebibliography}
\end{document}